\documentclass[12pt,a4paper]{article}
\usepackage[utf8]{inputenc}
\usepackage[T1]{fontenc}
\usepackage{amsmath,amssymb,amsfonts}
\usepackage{graphicx}
\usepackage{booktabs}
\usepackage{cite}
\usepackage{geometry}
\title{Constructing Pareto Sets Using Noether's Second Theorem}
\author{E.V. Nikulchev\\
MIREA -- Russian Technological University, Moscow, Russia}
\date{\today}

\begin{document}

\maketitle

\begin{abstract}
This paper presents an analytical approach to constructing the Pareto set in multiobjective variational control problems, based on Noether's second theorem. A fundamental connection is established between the gauge symmetries of the dynamical system describing the plant and the structure of the Pareto set. It is shown that the Pareto front can be interpreted as a conservation law arising from the invariance of the system and the quality criteria with respect to a symmetry group. A gauge regularization concept is proposed, which allows eliminating variables that do not affect the criteria without changing the Pareto set, thereby reducing the dimension of the solution space. Numerical examples are provided for the linear-quadratic regulator with two conflicting criteria and for the standard ZDT1 test problem.
\end{abstract}

\section{Introduction}

Many practical control and design problems involve simultaneous optimization of several conflicting criteria. In such cases, the solution is not a single point but a set of compromise solutions -- the Pareto set (Pareto front)~\cite{Taha2020}. Despite the large number of numerical methods for finding the Pareto front (evolutionary algorithms, weighted sum method, gradient approaches)~\cite{Li2009}, their efficiency is often limited by the high dimensionality of the solution space, constraints, and symmetries of the problem.

In this paper, we propose using Noether's second theorem to improve the efficiency of Pareto front search. At first glance, Noether's theorem belongs to mathematical physics, while the Pareto set belongs to control theory and decision making. However, both concepts express a linear dependence of gradients (variational derivatives) at an extremum point. In multiobjective optimization, the necessary condition for Pareto optimality is the existence of nonnegative Lagrange multipliers for which the weighted sum of gradients vanishes~\cite{Miettinen1999,Ehrgott2005}. In Noether's second theorem, the identity linear dependence of variational derivatives arises from gauge invariance and is expressed in the form of Bianchi identities~\cite{Olver1998}. The main ideas of this paper were developed in the author's PhD thesis in 2000~\cite{Nikulchev2000}, but new computational tools and modern algorithms require updating and rethinking the results obtained.

The connection between Noether's theorems and optimization problems has been studied in only a few works. An obvious connection between Noether's second theorem and the single-objective optimization problem at the analytical level was investigated in~\cite{Torres2003}, where D.F.M. Torres gives explicit expressions for Noether currents under gauge symmetries. In~\cite{Malinowska2009}, a multiobjective problem was considered for discrete models in the sense of~\cite{Censor1977}, where a ``scalarization of multiobjective optimization'' was introduced without using Noether's theorems. These studies did not lead to widespread use because they did not consider realistic control problems, lacked computational procedures and algorithms, and remained purely theoretical exercises with specialized mathematical apparatus. Thus, Noether's second theorem in the multiobjective context remains an unexplored area~\cite{Nikulchev2005}.

The main idea is as follows. The plant is described by a system of differential equations possessing symmetries. Several criteria are introduced, not abstractly, but specifically related to the plant in the form of Lagrangians. Thus, each criterion is invariant with respect to the same group. The main result is that the Pareto set is, in a certain sense, a conservation law: a compromise where improvement in one criterion inevitably leads to deterioration of the others. Constructing conservation laws is the essence of Noether's second theorem, according to which the variational derivatives of each criterion with respect to gauge variables are identically zero. This means that gauge variables do not affect the criteria values and can be either eliminated from the search space (dimension reduction) or fixed using regularization without changing the Pareto set. It is important to emphasize that Noether's theorem applies to each criterion individually and to the weighted sum. The proposed approach, although visually similar to linear convolution, derives the weights analytically from the symmetries of the problem rather than arbitrarily. Thus, an analytical method for constructing the Pareto front is obtained, which corresponds to a computationally reliable procedure that can be applied to modern Pareto front construction algorithms, significantly expanding their capabilities.

The rest of the paper is organized as follows: Section~2 presents the basic concepts and explains the essence of the results; Section~3 provides the mathematical formulation of the problem, the main theorem and its consequences for numerical algorithms; Section~4 contains the solution of a classical control problem with conflicting objectives; Section~5 describes the search for the Pareto set for a well-known test problem and compares the results with published accuracy data; Section~6 presents the conclusions.

\section{Main Result: Differential-Geometric Interpretation}

Before constructing formal models, we first consider several geometric concepts~\cite{Hall2013}. They do not require deep knowledge of differential geometry but help to see the unified structure of the problem and the connection between symmetries of the plant and the structure of the Pareto front.

Let the state of the system be described by variables $x = (x^1,\ldots,x^n)$ and control $u = (u^1,\ldots,u^m)$. The set of all possible values $(x,u)$ at each time instant forms a so-called fiber bundle -- a space where over each point $x$ there is a fiber of possible controls $u$. In other words, this is the set of all admissible pairs ``state--control''.

Quality functionals depend not only on the current values of $x,u$ but also on their time derivatives (velocities, accelerations, etc.). To describe all such dependencies uniformly, one introduces the jet space $J^r$ of order $r$. The coordinates in this space are the variables $t,x,u$ and, in general, derivatives of $x,u$ up to order $r$ inclusive. For example, the first-order jet space $(r=1)$ has coordinates $(t,x,u,\dot{x},\dot{u})$. The Lagrangian $L_i$ is a function on the jet space, i.e., it depends on $t,x,u$ and their derivatives. Such dependence arises in most applied problems: quality criteria are expressed as integrals of functions containing $x,\dot{x},u$.

The variational derivative $E_i(L_i)$ is a generalization of the gradient to the case of functions depending on derivatives. It shows how the functional $L_i$ changes under a small variation of the trajectory $x(t)$ taking into account the constraints imposed by the system of differential equations. Geometrically, the Euler--Lagrange equations $E_i(L_i)=0$ express the stationarity condition of the functional on the space of admissible trajectories.

Symmetries of the system are transformations that map trajectories to trajectories. They act on the jet space naturally: not only the variables themselves but also their derivatives are transformed (according to the rules of differentiation). The group of such transformations is the symmetry group $G$. The group $G$ is a Lie group, where the elements are transformations, the operation is the Lie bracket (commutator, i.e., sequential application of transformations), and the identity element is the identity transformation. If the Lagrangians $L_i$ are invariant under $G$ (up to a total derivative), then by Noether's theorem conservation laws arise.

In classical mechanics, the conservation of energy is a consequence of time homogeneity, and the conservation of momentum is a consequence of space homogeneity. These laws reflect the physical nature. Similarly, in multiobjective optimization, the Pareto set can be viewed as the result of a numerical compromise. Then, if all criteria are invariant under some transformation group, the Pareto front ceases to be random -- it becomes an invariant of this group.

Since all criteria refer to the same plant and are defined by the same dynamical system, the symmetry group of this system (the set of transformations preserving the equations) generates common symmetries for all criteria built on its solutions. In other words, the criteria may be different functions of state and control, but they are defined on the same solution space of a single system, and the symmetry group of this system acts on all criteria simultaneously.

For example, consider two conflicting criteria $f_1$ and $f_2$ defined on the solution space of one system. At a regular point of the Pareto front, the Karush--Kuhn--Tucker condition holds:
\[
\lambda_1 \nabla f_1 + \lambda_2 \nabla f_2 = 0, \quad \lambda_1,\lambda_2 > 0.
\]
This means that the gradients of the criteria are linearly dependent and point in opposite directions -- improvement in one criterion inevitably leads to deterioration of the other. This is a ``conservation law'' in the criterion space: the sum of weighted gradients is zero.

According to Noether's second theorem~\cite{Olver1998}, if functionals defined on the solutions of a system are invariant under a gauge symmetry group of that system, then their variational derivatives are related by Bianchi identities. These identities impose restrictions on the weights $\lambda_i$, allowing them to be eliminated and yielding a pure front equation -- an invariant independent of the choice of weights. Thus, the shape of the front is predetermined by the symmetries of the system and can be obtained analytically without solving any optimization problem.

Noether's second theorem provides a systematic way to find such invariants: gauge symmetries lead to identities relating variational derivatives, which in turn fix the admissible weights and lead to an explicit front equation.

The obtained relation looks similar to the classical weighted sum method with linear convolution $F = \xi_1 f_1 + \xi_2 f_2$. This similarity is deceptive: the weighted sum method requires solving optimization problems, does not provide an analytical compromise structure, and the weights $\xi_i$ are chosen empirically without physical justification.

The proposed approach is fundamentally different: the coefficients are derived from Noether's theorem and the symmetries of the system; regularization is applied to each criterion so that it does not change the Pareto set; dimension reduction (eliminating gauge variables) gives an analytical equation of the front or its parametrization; we obtain a single formula (invariant) that completely describes the front.

\section{Main Result: Mathematical Formalization}

Consider the plant described by a system of differential equations in normal form:
\[
\dot{x} = f(x,u), \quad x \in \mathbb{R}^n, \quad u \in \mathbb{R}^m, \tag{1}
\]
where $x$ is the state vector, $u$ is the control vector. Admissible pairs $(x(t),u(t))$ satisfy (1) on some time interval.

Let $G$ be a Lie group of transformations acting on the space $(t,x,u)$:
\[
\bar{t} = \varphi(t,x,u;\varepsilon), \quad \bar{x} = \psi(t,x,u;\varepsilon), \quad \bar{u} = \chi(t,x,u;\varepsilon),
\]
where $\varepsilon = (\varepsilon^1,\ldots,\varepsilon^r)$ are group parameters (in the gauge case, arbitrary functions of time). The group $G$ is called a \emph{symmetry group} of system (1) if it maps solutions to solutions, i.e., for any trajectory $(x(t),u(t))$ satisfying (1), the transformed trajectory $(\bar{x}(\bar{t}),\bar{u}(\bar{t}))$ also satisfies (1).

In infinitesimal form, the group generator is
\[
X = \tau(t,x,u)\frac{\partial}{\partial t} + \sum_{i=1}^n \xi^i(t,x,u)\frac{\partial}{\partial x^i} + \sum_{j=1}^m \eta^j(t,x,u)\frac{\partial}{\partial u^j}.
\]
The prolongation of $X$ to derivatives (including $\dot{x}$) is defined by standard Lie--Bäcklund formulas. The invariance condition of system (1) is written as
\[
X^{(1)}[\dot{x} - f(x,u)] \big|_{\dot{x}=f(x,u)} = 0. \tag{2}
\]
If the group $G$ depends on arbitrary functions of time (gauge group), the generator contains derivatives of these functions, and equation (2) leads to differential identities -- Bianchi identities -- relating the variational derivatives of functionals defined on solutions of the system.

Let $k$ functionals be defined on the solution space of system (1):
\[
J_i[u] = \int_\Omega L_i(x,u,u_{(1)},\ldots,u_{(r)})\, dt, \quad i=1,\ldots,k. \tag{3}
\]
The variational derivative (Euler--Lagrange equations) is
\[
E_i(L_i) = \frac{\delta J_i}{\delta u} = \sum_{|\alpha|\ge 0} (-1)^{|\alpha|} D_\alpha \frac{\partial L_i}{\partial u_\alpha}.
\]

\textbf{Definition~\cite{Malinowska2009}.} A trajectory $(x_0(t),u_0(t))$ is called weakly efficient (weakly Pareto optimal) if there is no variation $\delta u$ such that $J_i[u_0+\delta u] < J_i[u_0]$ for all $i$.

\textbf{Necessary condition for Pareto optimality} (generalized Karush--Kuhn--Tucker condition): if $u_0$ is weakly Pareto optimal, then there exist numbers $\lambda_i \ge 0$, not all zero, such that
\[
\sum_{i=1}^k \lambda_i E_i(L_i)(u_0) = 0.
\]
This result is given in~\cite{Miettinen1999,Ehrgott2005} for the finite-dimensional case; its generalization to infinite-dimensional variational problems is discussed in~\cite{Gorbunov2016}.

Assume that each functional (3) is invariant under the gauge group $G$ (a symmetry group of system (1)) depending on arbitrary functions $\varepsilon^a(x)$. Noether's second theorem states that there exist identities among the variational derivatives -- \emph{Bianchi identities}~\cite{Olver1998,Kosmann2011}:
\[
\sum_{i=1}^k \mu_i^{(a)}(u) E_i(L_i) \equiv 0, \quad a=1,\ldots,p, \tag{4}
\]
where $\mu_i^{(a)}$ are linear differential operators. In particular, if the group acts only on the gauge variables $z$ without affecting the effective variables, then the variational derivatives with respect to $z$ are identically zero for each criterion:
\[
\frac{\delta J_i}{\delta z} \equiv 0 \quad \forall i.
\]
This means that each criterion does not depend on the gauge variables.

\textbf{Theorem (Pareto Set as a Conservation Law).} Let a system of differential equations (1) describing the plant be invariant under a gauge Lie group $G$. Let functionals $J_i$, $i=1,\ldots,k$, each invariant under $G$, be defined on the solution space of this system. Then:

1. For any set $\lambda_i \ge 0$, every solution of the weighted Euler--Lagrange equation
\[
\sum_{i=1}^k \lambda_i E_i(L_i) = 0
\]
satisfies the necessary condition for Pareto optimality; the regular part of the Pareto set is contained in the union of such solutions over all $\lambda_i$.

2. On Pareto-optimal sections, combinations of Noether currents corresponding to symmetries of the system are conserved.

3. If the group $G$ is gauge (depends on arbitrary functions), then the Bianchi identities (4) impose restrictions on the weights $\lambda_i$, allowing them to be eliminated and yielding a pure front equation $\Phi(J_1,\ldots,J_k)=0$ independent of the weights.

\textbf{Proof.} Since all functionals are defined on the same solution space, invariant under $G$, the weighted sum $J_\lambda = \sum \lambda_i J_i$ is also defined on those solutions and is invariant. Substituting the Pareto condition into the Bianchi identities (4) yields a system of restrictions on the weights $\lambda_i$. The remaining statements follow from standard results of Noether's second theorem and linearity of the weighted sum. Details are given in~\cite{Olver1998}.

\textbf{Corollary for numerical algorithms.} From the invariance of the system and each criterion with respect to the gauge group $G$, a way to improve optimization algorithms follows~\cite{Beketov2025}. Gauge variables can be \emph{fixed by penalty}:
\[
F(x,z) = f_1(x) + \lambda f_2(x) + \mu \sum_{j=1}^{n_{\text{gauge}}} z_j^2,
\]
where $\mu > 0$ is the penalty parameter. According to Noether's theorem, such a penalty does not change the Pareto set. The parameter $\mu$ can be chosen arbitrarily (e.g., $\mu=1$), since its role is gauge fixing, not fitting to criteria. The penalty does not change the Pareto set because it vanishes at optimal solutions.

\section{Example: Linear-Quadratic Regulator with Conflicting Criteria}

Consider a linear system with $n$ states and one control:
\[
\dot{x} = A x + b u, \quad x(0) = x_0,
\]
where $x \in \mathbb{R}^n$, $u \in \mathbb{R}$, $A \in \mathbb{R}^{n\times n}$, $b \in \mathbb{R}^n$.

Instead of the traditional criterion as a sum, we consider two quadratic criteria:
\[
J_1 = \int_0^\infty x^T Q x \, dt, \quad J_2 = \int_0^\infty r u^2 \, dt,
\]
where $Q \succeq 0$ is the state weighting matrix, $r>0$ is the control weight. The first criterion minimizes the integral error (deviation of state from zero), the second minimizes the control energy. These goals conflict: reducing error requires increasing control, which raises $J_2$, and vice versa.

The system admits a one-parameter scaling transformation:
\[
t \rightarrow \alpha t, \quad x(t) \rightarrow x(\alpha t), \quad u(t) \rightarrow \alpha^{-1} u(\alpha t), \quad \alpha > 0.
\]
In infinitesimal form for $\alpha = 1+\varepsilon$, $\varepsilon \ll 1$, the variations are
\[
\delta t = \varepsilon t, \quad \delta x = 0, \quad \delta u = -\varepsilon u.
\]
The corresponding infinitesimal generator acting on functions $(t,x,u)$ is
\[
X = t \frac{\partial}{\partial t} - u \frac{\partial}{\partial u}.
\]

Consider the parametric family of problems depending on the additional variable $\lambda$; Noether's second theorem will give a conservation law. We derive the Bianchi identity.

Consider the functional
\[
J_\lambda = J_1 + \lambda J_2 = \int_0^\infty \left( x^T Q x + \lambda r u^2 \right) dt.
\]
Consider the infinitesimal transformation acting on the parameter $\lambda$: $\lambda \rightarrow \lambda + \varepsilon$, $\varepsilon \ll 1$. In the extended space of variables $(t,x,u,\lambda)$, this transformation has the form:
\[
\delta \lambda = \varepsilon.
\]
The corresponding generator is
\[
X_\lambda = \frac{\partial}{\partial \lambda}.
\]
This generator commutes with $X$.

Compute the variation of $J_\lambda$ under $\lambda \rightarrow \lambda+\varepsilon$. Since $x$ and $u$ are not varied, only $\lambda$ changes, we have
\[
\delta J_\lambda = \varepsilon \int_0^\infty r u^2 \, dt = \varepsilon J_2.
\]
On the other hand, if we consider $J_\lambda$ as a function of the parameter $\lambda$ on optimal trajectories, its total derivative with respect to $\lambda$ equals $J_2(\lambda)$, where $J_2(\lambda)$ is the value of the second criterion at the optimal solution for that $\lambda$. The Euler--Lagrange equation for $L_\lambda$ is
\[
-\frac{d}{dt}\left( \frac{\partial L_\lambda}{\partial \dot{x}} \right) + \frac{\partial L_\lambda}{\partial x} = 0,
\]
and its solution $x_\lambda(t)$ depends on $\lambda$. Differentiating $J_\lambda$ with respect to $\lambda$ and using the Euler--Lagrange equation gives
\[
\frac{d}{d\lambda} J_\lambda = \int_0^\infty \frac{\partial L_\lambda}{\partial \lambda} dt = \int_0^\infty r u_\lambda^2 dt = J_2(\lambda).
\]
This relation is the Bianchi identity for the parametric problem.

Since $J_\lambda = J_1 + \lambda J_2$, differentiating with respect to $\lambda$ and using the above, we get
\[
\frac{dJ_1}{d\lambda} + J_2 + \lambda \frac{dJ_2}{d\lambda} = J_2,
\]
which yields the differential conservation law relating changes in criteria when varying the weight $\lambda$:
\[
\frac{dJ_1}{d\lambda} + \lambda \frac{dJ_2}{d\lambda} = 0. \tag{5}
\]

\textbf{Computational Example 1.} Consider the scalar case $n=1$, $r=1$.

For $n=1$, $A=0$, $b=1$, $Q=1$, $r=1$, the system is $\dot{x}=u$, criteria:
\[
J_1 = \int_0^\infty x^2 dt, \quad J_2 = \int_0^\infty u^2 dt.
\]
In this case, the conservation law integrates:
\[
\frac{d}{d\lambda}(J_1 - \lambda J_2) = 0 \quad \Rightarrow \quad J_1 - \lambda J_2 = \text{const}.
\]
From the analytical solution of LQR using the Riccati equation for this system: $J_1 = x_0^2/(2\sqrt{\lambda})$, $J_2 = x_0^2/(2\lambda^{3/2})$, we obtain $J_1 - \lambda J_2 = 0$. Hence the constant is zero, and the conservation law becomes
\[
J_1 = \lambda J_2. \tag{6}
\]
This is an algebraic invariant following from the symmetry of the problem.

Figure~1 shows the analytical dependence and 50 Pareto set points for the scalar case for $\lambda$ ranging from 0.01 to 10 on a logarithmic scale (for uniform distribution on the hyperbola).

\begin{figure}[htb]
\centering
\includegraphics[width=0.8\textwidth]{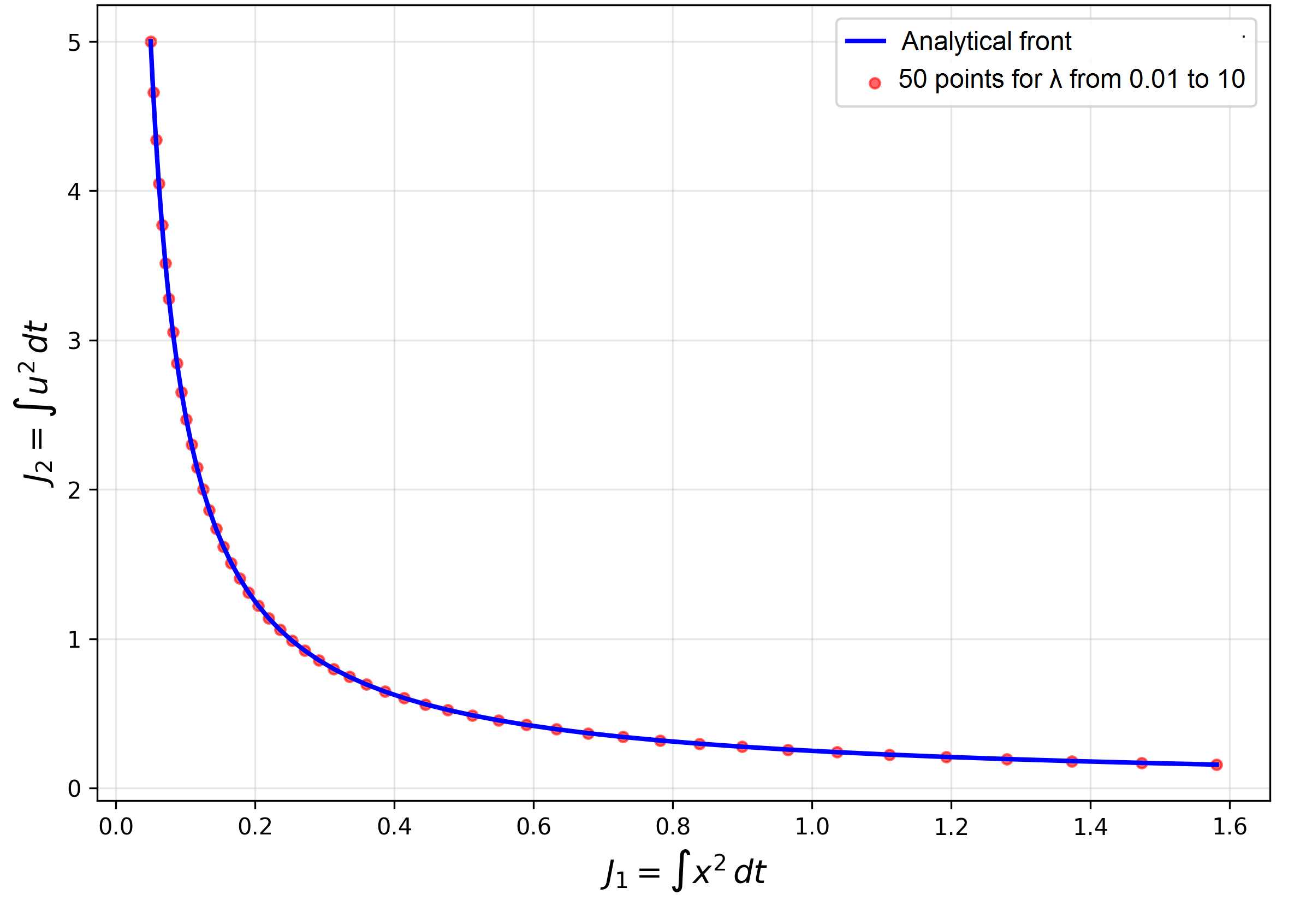}
\caption{Analytical front for the scalar LQR problem obtained from the conservation law.}
\label{fig:lqr1}
\end{figure}

\textbf{Example 2.} State space dimension 3.

For the system with $n=3$ matrices
\[
A = \begin{bmatrix}
-1 & 0.5 & 0.2 \\
0.3 & -2 & 0.1 \\
0 & 0.4 & -1.5
\end{bmatrix}, \quad
b = \begin{bmatrix}
1 \\ 0.5 \\ 0
\end{bmatrix}, \quad
x_0 = \begin{bmatrix}
1 \\ 0 \\ 0
\end{bmatrix}, \quad
Q = I_3, \quad r=1,
\]
the Pareto front was constructed for $\lambda \in [0.01, 10]$ using the Riccati equation. For each $\lambda$, the equation
\[
A^T P + P A - \frac{1}{\lambda r} P b b^T P + Q = 0
\]
was solved, $J_\lambda(\lambda) = x_0^T P(\lambda) x_0$ was computed, and by central differences with step $\varepsilon=10^{-4}$, $J_1(\lambda)$ and $J_2(\lambda)$ were found via the conservation law:
\[
J_2(\lambda) = \frac{d}{d\lambda} J_\lambda(\lambda), \quad J_1(\lambda) = J_\lambda(\lambda) - \lambda J_2(\lambda).
\]
Figure~2 shows the obtained Pareto front.

\begin{figure}[htb]
\centering
\includegraphics[width=0.8\textwidth]{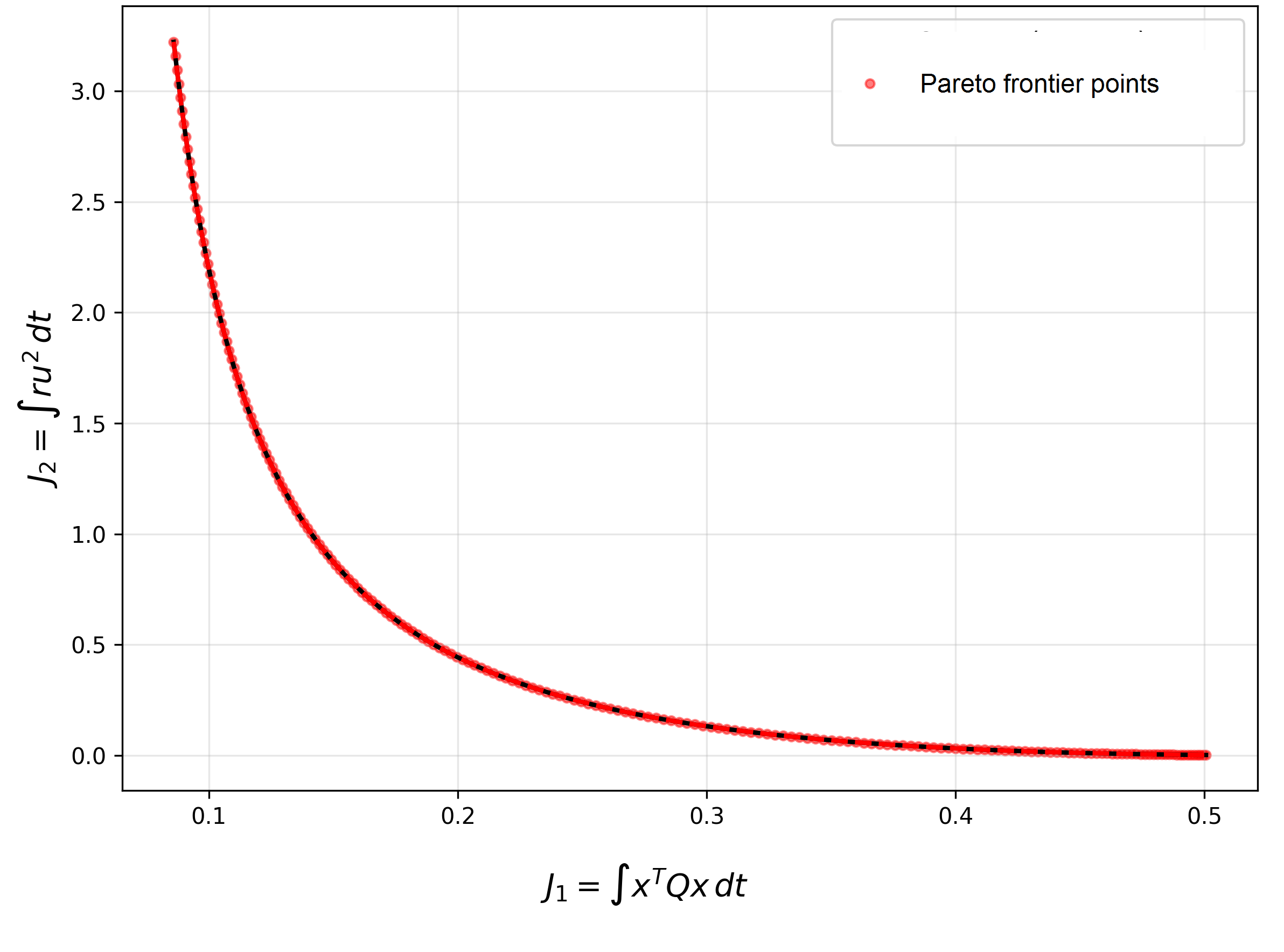}
\caption{Pareto front for Example 2.}
\label{fig:lqr2}
\end{figure}

In the classical approach without using the conservation law, for each $\lambda$ one must numerically solve an optimal control problem, e.g., by time discretization and nonlinear programming. For comparison, computations were performed in Colab using SciPy Optimize. For the system of dimension $n=3$ with horizon $T=20$ and discretization step $\Delta t = 0.01$, solving one problem for four values of $\lambda$ took about 1500 seconds, and constructing 100 front points required nearly one and a half hours. The correctness of the front depends on the discretization step. Figure~3 shows 4 points found with discretization steps 0.1 and 0.5 compared with the front obtained from the conservation law. The method based on the conservation law took only 0.2 seconds.

\begin{figure}[htb]
\centering
\begin{minipage}[b]{0.45\textwidth}
\centering
\includegraphics[width=\textwidth]{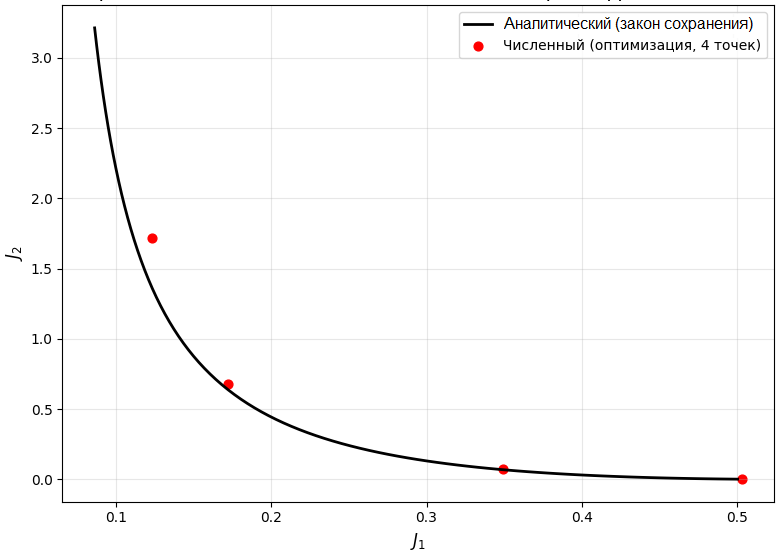}
\small (a) Discretization step 0.1
\end{minipage}
\hfill
\begin{minipage}[b]{0.45\textwidth}
\centering
\includegraphics[width=\textwidth]{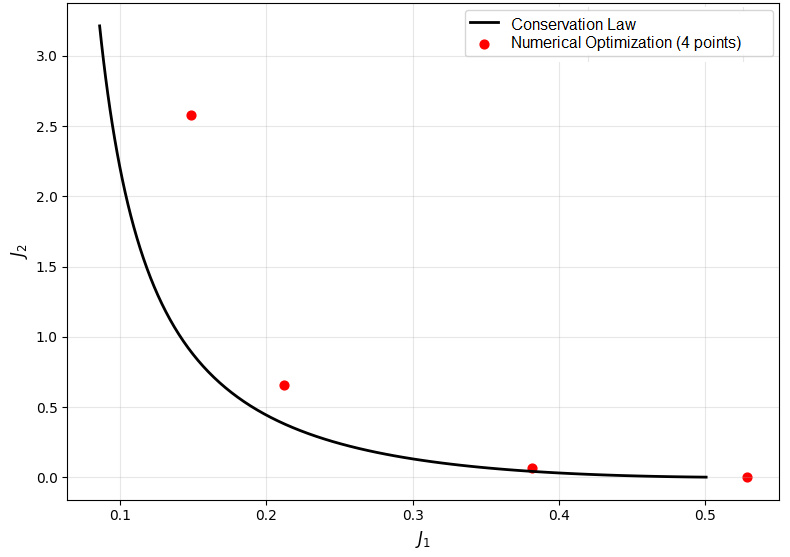}
\small (b) Discretization step 0.01
\end{minipage}
\caption{Comparison of the front from the conservation law (solid line) and numerically computed points.}
\label{fig:lqr3}
\end{figure}

Thus, the conservation law arising from parametric symmetry of the problem allows constructing the Pareto front analytically, without repeatedly solving the optimization problem, giving a huge gain in speed and accuracy.

The examples illustrate the application of the proven general theorem on the connection between symmetries and the structure of the Pareto front and show the practical advantage of the analytical method over numerical ones.

\section{Test Problem: ZDT1}

To evaluate the quality of Pareto front search algorithms, a set of standard test problems is used, where the analytical form of the set is known. This allows comparing the obtained results with published ones, for example, in terms of accuracy.

Consider the ZDT1 problem~\cite{Zitzler2000}:
\[
f_1(x) = x_1, \quad f_2(x) = g(x) \left( 1 - \sqrt{x_1 / g(x)} \right),
\]
where
\[
g(x) = 1 + \frac{9}{n-1} \sum_{i=2}^{n} x_i, \quad x_i \in [0,1].
\]
The theoretical front is $f_2 = 1 - \sqrt{f_1}$.

Note that this problem has no physical meaning; it is a test for algorithm evaluation.

Introduce the scalarized functional
\[
F_\lambda(x) = f_1(x) + \lambda f_2(x), \quad \lambda > 0.
\]
Consider a one-parameter group of transformations acting on the variables $(x,\lambda)$:
\[
\bar{x}_1 = \frac{\lambda}{\mu} x_1, \quad \bar{x}_i = x_i \ (i \ge 2), \quad \bar{\lambda} = \mu \lambda,
\]
where $\mu > 0$ is the group parameter. Direct substitution verifies that
\[
F_{\bar{\lambda}}(\bar{x}) = F_\lambda(x),
\]
i.e., the functional is invariant under this transformation. In infinitesimal form for $\mu = 1+\varepsilon$, $\varepsilon \ll 1$, the generator is
\[
X = \lambda \frac{\partial}{\partial \lambda} + x_1 \frac{\partial}{\partial x_1}.
\]
The prolongation of the generator to derivatives (in this case, derivatives with respect to $\lambda$) gives
\[
X^{(1)} = X + \dot{x}_1 \frac{\partial}{\partial \dot{x}_1} + \dot{\lambda} \frac{\partial}{\partial \dot{\lambda}},
\]
where $\dot{x}_1 = dx_1/d\lambda$, $\dot{\lambda}=1$.

The invariance condition of $F_\lambda$ with respect to $X$ is written as $X^{(1)} F_\lambda = 0$ on solutions of the Euler--Lagrange equation. Expanding this expression yields the Bianchi identity:
\[
\frac{\partial F_\lambda}{\partial \lambda} + \frac{\partial F_\lambda}{\partial x_1} \frac{dx_1}{d\lambda} = 0.
\]
On the optimal trajectory (where $\partial F_\lambda/\partial x_1 = 0$), the conservation law follows:
\[
f_2(x) + \lambda \frac{df_2}{d\lambda} + \frac{df_1}{d\lambda} = 0. \tag{7}
\]

In the numerical implementation, to construct an approximation of the front, random initial points were used, and gauge variables $z$ were introduced with a quadratic penalty $\mu \|z\|^2$ (here $\mu = 10$). A total of 500 front points were constructed. The average execution time over 30 runs was 140 seconds. Results are shown in Figure~4.

\begin{figure}[htb]
\centering
\includegraphics[width=0.8\textwidth]{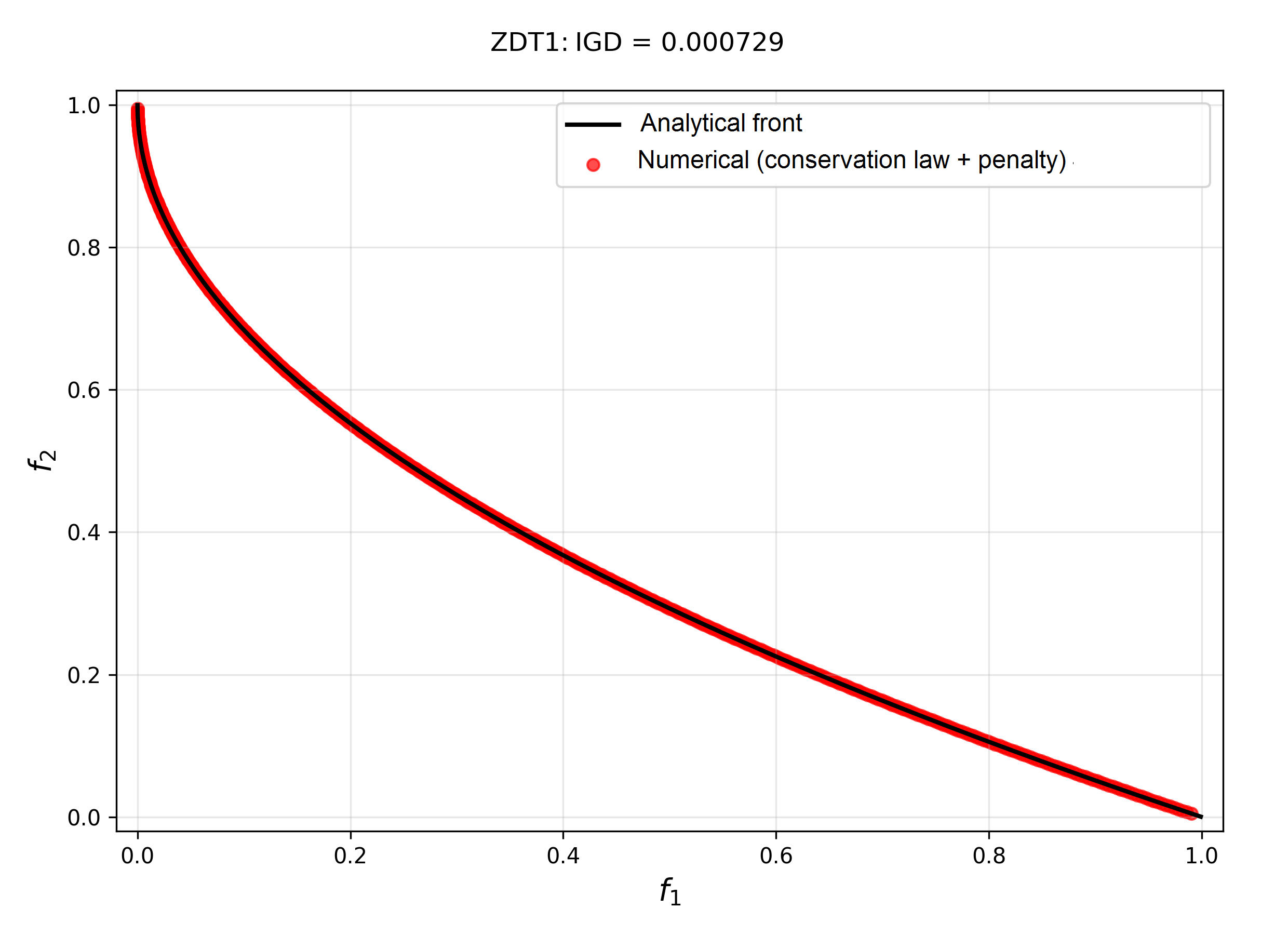}
\caption{Comparison of the front found via the conservation law (points) with the analytical front (solid line) for ZDT1.}
\label{fig:zdt1}
\end{figure}

To evaluate the accuracy of the Pareto front approximation, the Inverted Generational Distance (IGD) metric is used:
\[
IGD(P^*,P) = \frac{1}{|P^*|} \sum_{p^* \in P^*} \min_{p \in P} \|p^* - p\|_2,
\]
where $P^*$ is the set of reference front points and $P$ is the set of points found by the algorithm, and $\|\cdot\|_2$ is the Euclidean norm in the criterion space. Although it has been shown in~\cite{Rodriguez2025} that the choice of reference set critically affects indicator values and using a biased reference set can lead to incorrect quality assessments, this metric is widely used in the literature.

A comparison of IGD with published results is given in Table~\ref{tab:igd}. The results show that the proposed approach achieves IGD = $0.729\times 10^{-3}$, which is the best value for ZDT1 among the reported publications. This confirms that regularization based on gauge variables, justified by Noether's second theorem, allows achieving high accuracy in Pareto front approximation on standard test problems.

\begin{table}[htb]
\centering
\caption{Comparison of algorithm accuracy (IGD) for ZDT1.}
\label{tab:igd}
\begin{tabular}{lcc}
\toprule
\textbf{Approach} & \textbf{Algorithm} & \textbf{IGD} \\
\midrule
Noether's theorem & Conservation law & $\mathbf{0.729 \times 10^{-3}}$ \\
Particle swarm & MMOPSO~\cite{Chen2025} & $2.446\times 10^{-3}$ \\
& MOPSOCD~\cite{Luo2025} & $2.804\times 10^{-3}$ \\
& ASDMOPSO~\cite{Ye2023} & $2.918\times 10^{-3}$ \\
& NMPSO~\cite{Chen2025} & $25.909\times 10^{-3}$ \\
Evolutionary decomposition & MOEADCMA~\cite{Ye2023} & $3.899\times 10^{-3}$ \\
& MOEA/D~\cite{Kalita2024} & $32.328\times 10^{-3}$ \\
& MOEDO~\cite{Kalita2024} & $5.753\times 10^{-3}$ \\
NSGA-II & NSGA-II~\cite{Ye2023} & $4.653\times 10^{-3}$ \\
& NSGA-II~\cite{Kalita2024} & $4.111\times 10^{-3}$ \\
Other & NSCSO~\cite{Huang2024} & $4.93\times 10^{-3}$ \\
& LMEA~\cite{Huang2024} & $5.33\times 10^{-3}$ \\
\bottomrule
\end{tabular}
\end{table}

\section{Conclusion}

In this work, an approach to constructing the Pareto set in multiobjective variational control problems based on Noether's second theorem is substantiated. A fundamental connection between gauge symmetries of the dynamical system describing the plant and the structure of the Pareto set is established: the Pareto front can be interpreted as a conservation law arising from the invariance of the system and criteria with respect to a symmetry group.

A theorem is proved stating that if the system of differential equations and each quality criterion are invariant under a gauge Lie group, then the Pareto set is an invariant of this group, and its shape can be described analytically through Bianchi identities. Differential conservation laws relating changes in criteria under variation of weight coefficients are obtained.

It is shown that gauge variables that do not affect the criteria values can be eliminated from the search space or fixed with a quadratic penalty without changing the Pareto set.

The proposed method is tested on the classical linear-quadratic regulator problem with two conflicting criteria and on the standard ZDT1 test problem. The results show that the analytical approach based on conservation laws allows constructing the Pareto front hundreds of times faster than direct numerical solution of optimal control problems, with accuracy determined only by the discretization step when computing derivatives. For ZDT1, the developed algorithm achieved IGD = $0.729\times 10^{-3}$, which surpasses published results for modern algorithms, including MMOPSO ($2.446\times 10^{-3}$), MOEADCMA ($3.899\times 10^{-3}$), and standard NSGA-II ($4.653\times 10^{-3}$).

The obtained results may open a new direction in multiobjective optimization, allowing not only to improve existing numerical algorithms but also to obtain analytical solutions for a wide class of problems with symmetries. The proposed approach can be extended to more complex control systems, including nonlinear and time-varying plants. Future research includes the development of specialized hybrid algorithms combining analytical conservation laws with evolutionary methods, as well as application of the proposed approach to real engineering problems.

\end{document}